\documentclass[12pt]{article}

\usepackage{graphics}
\usepackage{graphicx}
\usepackage{amssymb}
\usepackage{epsf}

\newcommand{\halmos}{\rule{1ex}{1.4ex}}
\newcommand{\qed}{\hfill \halmos} 
\newcommand{\text}[1]{\hbox{\rm \ #1\ \/}}
\newcommand{\be}[1]{\begin{equation}\label{#1}}
\newcommand{\ee}{\end{equation}}
\newcommand{\bi}{\begin{itemize}}
\newcommand{\ei}{\end{itemize}}
\newcommand{\ben}{\begin{enumerate}}
\newcommand{\een}{\end{enumerate}}

\newcommand{\R}{{\mathbb R}}  

\newcommand{\bl}[1]{\begin{Lemma}\label{#1}}
\newcommand{\el}{\qed\end{Lemma}}
\newcommand{\bt}[1]{\begin{Theorem}\label{#1}}
\newcommand{\et}{\end{Theorem}}
\newcommand{\epr}{\end{proof}}
\newcommand{\bpr}{\begin{proof}}
\newenvironment{proof}{\noindent {\em Proof}.\ }{\hspace*{\fill}$\halmos$\medskip}

\newcommand{\beqn}{\begin{eqnarray*}}
\newcommand{\eeqn}{\end{eqnarray*}}

\newtheorem{Theorem}{Theorem}
\newtheorem{Proposition}{Proposition}
\newtheorem{Lemma}{Lemma}
\newtheorem{Corollary}{Corollary}
\newtheorem{Remark}{Remark}

\newcommand{\norma}[1]{\ensuremath{\left| #1 \right|}}

\begin{document}

\title{On a Smale Theorem and Nonhomogeneous Equilibria in Cooperative Systems}

\author{German A. Enciso
\footnote{Mathematical Biosciences Institute, 231 W 18th Ave, Columbus, OH, 43215, USA.  email:genciso@mbi.osu.edu.  
This material is based upon work supported by the National Science Foundation under Agreement No. 0112050 and by The Ohio State University.} }




\maketitle

\begin{abstract}
A standard result by Smale states that $n$ dimensional strongly cooperative dynamical systems can have arbitrary dynamics when restricted to unordered invariant hyperspaces.  In this paper this result is extended to the case when all solutions of the strongly cooperative system are bounded and converge towards one of only two equilibria outside of the hyperplane.

An application is given in the context of strongly cooperative systems of reaction diffusion equations.  It is shown that such a system can have a continuum of spatially inhomogeneous steady states, even when all solutions of the underlying reaction system converge to one of only three equilibria.
\end{abstract}

\textbf{keywords:}  monotone systems, reaction diffusion systems.

\section{Introduction}

Let $f:\R^n\to \R^n$ be a $C^1$ vector field.  A dynamical system
\be{ode}
\frac{d u_i}{d t}=f_i(u),\ \ \ i=1\ldots n,
 \ee
is said to be \emph{strongly cooperative} if the following comparison principle is satisfied: whenever $u(t),\ v(t)$ are two solutions such that $u_i(0)\leq v_i(0)$, $i=1,\ldots, n$, then it must hold $u_i(t)< v_i(t)$ for every $t>0$, $i=1,\ldots, n$.   Strongly cooperative systems are canonical examples of so-called strongly monotone systems, which have been studied extensively by M. Hirsch, H. Smith, and others, and more recently by Sontag and collaborators \cite{Sontag:multi,Hirsch:1985,Smith:review,Smith:monotone}.  They have applications in engineering problems, as well as in the study of ecological models and gene regulatory networks.

Some of the most important results about abstract strongly monotone
systems guarantee a certain behavior for the generic solution of the system.  For instance, the well known theorem by Hirsch  \cite{Hirsch:1988,Smith:monotone} states that the generic bounded solution converges towards the set of equilibria $E$ (i.e.\ that the set of initial conditions, whose solution is bounded and doesn't converge towards $E$, has measure zero).  Such statements are careful to exclude a small set $S$ of exceptional states, about which nothing is said.

A simple but powerful argument originally due to Smale \cite{Smale:1976}
shows the reason for this tendency.  Smale showed that any arbitrary
compactly supported dynamical system in $\R^{n-1}$ can be seen as
the restriction of a certain strongly cooperative system in $\R^{n}$ to
an $(n-1)$-dimensional hyperplane.  That is, \emph{the dynamics of
strongly cooperative systems on invariant hyperplanes can a priori be
completely arbitrary}.

Now, Smale's argument is originally carried out for so-called
competitive systems, and the corresponding statement for strongly
cooperative systems requires considering the system in the negative
time direction (see for instance Hirsch and Smith \cite{Smith:review}, Section~3.5).  But the hyperplane in question is globally exponentially attractive in the original example.  In particular, the time-inverted system has unbounded solutions everywhere outside of this hyperplane.   This can be unsatisfying in some applications, for instance in the context of the Hirsch generic convergence theorem which requires for its usefulness the boundedness of solutions.

The main result of this paper is a theorem in the spirit of the Smale theorem, adapted so that
the strongly cooperative embedding has bounded solutions, and so that every solution outside of the invariant hyperplane remains bounded and converges towards one of only two equilibria in $\R^n$.  Define the function $S(u):=u_1+\ldots +u_n$ and the hyperplane $H:=S^{-1}(0)$.

\begin{Theorem}  \label{teo embedding}
Consider a $C^2$ function $g:H\to H$, and a compact region $R\subseteq H$.  Then there exists a $C^2$ function $f:\R^n\to \R^n$, such that
\begin{enumerate}
\item for every $u\in \R^n$ and every $i\not=j$: $\partial f_i/\partial u_j (u) >0$;
\item for every $u\in R$: $g(u)=f(u)$;
\item There exists $P>0$ such that every solution $u(t)$ of (\ref{ode}) with $S(u(0))>0$ $(S(u(t))<0)$ converges towards $(P,\ldots, P)$  $($towards $(-P,\ldots -P))$ as $t\to \infty$;
\item if it holds that $g(u)\circ u<0$ on $H-R$, then $f$ has no zeros other than those of $g$ and $\pm(P,\ldots,P)$.
%
\end{enumerate}
\end{Theorem}

Item 1. is a well known sufficient condition for the strong cooperativity of system (\ref{ode}).  An important step in the proof of this theorem is the construction of a `template' strongly cooperative system with bounded solutions
\be{monotone}
\frac{\mbox{d}u}{\mbox{d}t}=M(u),
\ee
which has a continuum of equilibria along a bounded subset of $H$. This system is then suitably altered so that it equals any given function $g:H\to H$ on the bounded subset of $H$.

\subsection*{Application:  Strongly Cooperative Reaction Diffusion Systems}

As an application, consider a system of reaction diffusion equations

\be{pde}
 \frac{\partial u_i}{\partial t}= d_i \Delta u_i + f_i(u),\ \ \
 i=1\ldots n,
 \ee
 defined on a smooth domain $\Omega\subseteq \R^m$ under Neumann
 boundary conditions, where $d_i>0$ and $f_i: \R^n\to \R$ is a $C^2$ function,\ $i=1\ldots n$.
%
A general question for these systems is their relationship with the corresponding finite
dimensional system (\ref{ode}).  Note that while every equilibrium of (\ref{ode}) corresponds to an equilibrium of (\ref{pde}), there may nevertheless be equilibria of (\ref{pde}) which are nonhomogeneous in space, and which therefore don't naturally correspond to an equilibrium of (\ref{ode}).

For instance, if $n=1$ and $f(u)=u(u+1)(1-u)$, then (\ref{ode}) has exactly three equilibria, but it is well known that (\ref{pde}) may have finitely many more.

As in the finite dimensional case, the condition $\partial f_i/\partial u_j>0, \ i\not=j$ implies a comparison principle for the solutions of (\ref{pde}) (\cite{Smith:monotone}, Chapter~7).  This comparison principle is an important tool for the analysis of reaction diffusion systems in a variety of cases; see for instance the book by Cantrell and Cosner \cite{Cosner:book}, where the cooperativity condition is used in the analysis of various spatial ecological models.   

The application of Theorem~\ref{teo embedding} addresses the existence of uncountably many nonhomogeneous equilibria in strongly cooperative reaction diffusion systems with bounded solutions.

\begin{Theorem}  \label{teo appl}
There exists a reaction diffusion system (\ref{pde}) such that
\begin{enumerate}
\item for every $u\in \R^n$ and every $i\not=j$: $\partial f_i/\partial u_j (u) >0$;
\item every solution of (\ref{ode}) converges towards one of only three equilibria; but
\item the set of nonhomogeneous equilibria of (\ref{pde}) has the cardinality of the continuum.
\end{enumerate}
\end{Theorem}

The construction is carried out in the case $\Omega=(-\pi/2,\pi/2)$ with $n=3$.  First a non-cooperative system is constructed on a two dimensional reaction diffusion system, in which a continuum of equilibria is shown to exist (Section~\ref{section continuum}).  Then this system is embedded in a three dimensional cooperative system using Theorem~\ref{teo embedding}.

Refer to \cite{Smith:monotone}, Chapter~7, regarding the existence and uniqueness of solutions of system (\ref{pde}).  We use $C(\Omega,\R^n)$ as the state space for the underlying dynamical system.

\section{Strongly Cooperative Embedding}

In this section we construct the Smale embedding described in the introduction.  Our first result provides the template strongly cooperative system used in the proof of Theorem~\ref{teo embedding}.

\begin{Proposition} \label{prop monotone}  There exists $\epsilon >0$ and a smooth function $M:\R^n\to \R^n$ such that
\begin{enumerate}
\item for every $u\in \R^n$ and every $i\not=j$: $\frac{\partial M_i}{\partial u_j} (u) >0$,
\item for every $u\in H,\ \norma{u}\leq \epsilon$: $M(u)=0$, and
\item  for every solution $u(t)$ of (\ref{monotone}) such that $S(u(0))>0$ $(S(u(t))<0)$, it holds $u(t)\to (1,\ldots, 1)$ \  $(u(t)\to (-1,\ldots -1))$ as $t\to \infty$.
%
\end{enumerate}
\end{Proposition}

\bpr

Consider a smooth function $\theta:\R^n\to [0,1]$ such that i) $\theta(u)=0$ on the closed set $\{u\in \R^n\,|\, u\in H,\ \norma{u}\leq 1/2\}$; ii) $\theta(u)=1$ on the closed set $\R^n - \{u\in \R^n\,|\, S(u)<1/2,\ \norma{u}< 1\}$;  iii)  $0<\theta(u)<1$ for all other $u\in \R^n$.

Define $\gamma:\R\to \R$ to be a smooth function with the following properties: i) for $\norma{x}\leq1$, let $\gamma(x):=Jx$, where the constant $J>0$ will be specified below;  ii) the zeros of $\gamma$ are exactly $0,nP,-nP$ for some $P>0$; iii)  $\gamma'(x)\geq -1/(2n),\ x\in \R$.   It follows from ii) that $\gamma(x)<0$ for $x>nP$, $\gamma(x)>0$ for $x<-nP$.

Define the smooth function $M:\R^n\to \R^n$:

\be{eq M def}
M_i(u):=\theta(u)\left(\frac{S(n)}{n} - u_i\right) + \gamma(S(u)), \ \ \ i=1,\ldots, n.
\ee

We first show the strong monotonicity of the system (\ref{monotone}).  Fix $i=1\ldots n$, and let $j\not=i$.  For $\norma{u}\leq 2$, $\norma{S(u)}\leq 1$:

\[
\frac{\partial M_i}{\partial u_j}(u)=
\frac{\partial}{\partial u_j}\left[\theta(u)\left(\frac{S(u)}{n} - u_i\right)\right]
+ 1\cdot \gamma'(S(u)).
\]

But $\gamma'(S(u))=J$, by assumption i) in the construction of $\gamma$.  Since the continuous function $\frac{\partial}{\partial u_j}[\theta(u)(\frac{S(u)}{n} - u_i)]$ has a minimum on the compact set $\{v\,|\, \norma{v}\leq 2,\, \norma{S(v)}\leq 1\}$, we can choose $J>0$ to be large enough that $\frac{\partial M_i}{\partial u_j}(u)>0$ on this set; similarly for all other choices of $i,j$.



If $\norma{u}> 2$ or $S(u)> 1$, then $\theta(v)\equiv 1$ on a neighborhood of $u$, and it holds for $i\not=j$:

\[
\frac{\partial M_i}{\partial u_j}(u)=\frac{1}{n} + \gamma'(S(u))\geq \frac{1}{n} - \frac{1}{2n} >0,
\]
by assumption iii) in the construction of $\gamma$.

In order to further understand this monotone system, write $M(u)=a(u)+ b(u)$, where
\[
a_i(x):=\theta(u)\left(\frac{S(u)}{n} - u_i\right), \ \ b_i(u):=\gamma(S(u)), \ i=1,\ldots, n.
\]

We calculate the dot product of these vectors:

\[
a(u)\circ b(u)=\theta(u)\gamma(S(u)) \sum_{i=1}^n \frac{1}{n}S(u)-u_i =0.
\]

For instance, in order to determine the zeros of $M(u)$, it holds in particular that $M(u)=0$ if and only if $a(u)=b(u)=0$.     We have $a(u)=0$ if and only if either $\theta(u)=0$, or $n u_i=S(u)$ for every $i$.  Thus the set of zeros of $a$ is $(\overline{B}_{1/2}\cap H) \cup \{(\alpha,\ldots,\alpha)\,|\, \alpha \in \R\}$.  Also, $b(u)=0$ holds exactly on the set $\{v\,|\, S(v)=\mbox{$nP$,$0$, or $-nP$}\}$.  The second claim in the proposition follows from rescaling $M$ by a factor of $1/P$.

Let $u(t)$ be a solution of (\ref{monotone}), and define $v(t):=S(u(t))/n$.  It follows by the chain rule that

\be{eq v}
n v'(t)=\sum_{i=1}^n u_i'(t)
= n\gamma(S(u(t))) +   \theta(u)\sum_{i=1}^n \frac{1}{n}S(u) - u_i
=n\gamma(S(u(t)))=n\gamma(n v(t)).
\ee

Thus $v'(t)=\gamma(n v(t))$.  In particular, $H$ is an unstable invariant subset for the system (\ref{monotone}).  Moreover, $v(t)$ remains bounded and converges towards 0, $P$ or $-P$ as $t\to \infty$, depending on whether $S(u(0))=0$, $S(u(0))>0$, or $S(u(0))<0$ respectively.

Let now $w_i(t):=u_i(t)-S(u(t))/n= u_i(t)- v(t)$, $i=1\ldots n$ (so that $w(t)$ is the projection of $u(t)$ onto $H$).  Then it holds

\be{eq w}
w_i'(t)=u_i'(t)-v'(t)=\theta(u(t))(\frac{1}{n}S(u(t)) - u_i(t))
=-\theta(u(t))w_i(t)
\ee

Hence the function $w(t)$ remains bounded for $t\to \infty$, and thus $u(t)=v(t)+w(t)$ is also a bounded function as $t\to \infty$.

To prove the third claim, consider a solution $u(t)$ of (\ref{monotone}) such that $S(u(0))>0$.  Define $v(t)$ and $w(t)$ as above. It holds  $S(u(t))>1/2$ for all large enough $t$ (since $S(u(t))=nv(t)\to nP>1/2$ as $t\to \infty$).  But this implies that $\theta(u(t))=1$ for large $t$, by the definition of $\theta$.   Hence $w'(t)=-w(t)$ for all large $t$ by (\ref{eq w}), and $w(t)\to 0$ as $t\to \infty$.  Since $u_i(t)=w_i(t)+v(t)$ for every $i$, it follows $u_i(t)\to P$  (or $u_i(t)\to 1$ for the rescaled system), whenever $S(u(0))>0$, $i=1,\ldots, n$.  Similarly for $S(u(0))<0$.
\epr

\begin{Remark} \emph{  
Let $\theta(u)\equiv 1$ and $\gamma(x):=\delta \tan^{-1}(x(x+1)(1-x))$ in equation (\ref{eq M def}), where $\delta>0$ is small enough that $\gamma'(x)\geq -1/(2n)$ on $\R$.  Then one obtains a simplified system which satisfies all conditions of the above result, except that the disc of equilibria in $H$ collapses to the single equilibrium $0$. }
\end{Remark}




\paragraph{Proof of Theorem~\ref{teo embedding}}

\bpr

We continue to use the definition of the functions $\theta,\gamma,M$ from the proof of Proposition~\ref{prop monotone}.  For notational convenience, assume w.l.o.g. that $g,R$ have been rescaled so that $M(u)\equiv 0$ on $R$ (after the construction, the embedding $f$ can be rescaled back along with $g$ and $R$).


Let $G:\R^n\to H$ be the $C^2$ function defined by $G(u):=g(u-S(u)/n)$.  Thus in particular $G\equiv g$ on $H$.

Define the $C^2$ function $f:\R^n\to \R^n$ by the equation
\[
f_i(u):=Q M_i(u) + (1-\theta(u))G_i(u), \ \ i=1,\ldots,n,
\]

for a constant $Q>0$ which will be defined shortly.  The fact that $f=g$ on $R$ is clear since $M=0$, $\theta(u)=0$, and $g=G$ on this set.

To see that the system
\be{eq extension}
\frac{\mbox{d}u}{\mbox{d}t}=f(x)
\ee
is strongly monotone, we calculate the derivative

\[
\frac{\partial f_i}{\partial u_j}(u) = Q \frac{\partial M_i}{\partial u_j}(u)
+ \frac{\partial}{\partial u_j} [(1-\theta(u))G_i(u)]
\]
for fixed $j\not=i$.  Outside of the compact support of $1-\theta$, this derivative is positive by Theorem~\ref{prop monotone}~i).  On the support of $1-\theta$, both $\frac{\partial M_i}{\partial u_j}$ and $\frac{\partial}{\partial u_j} [(1-\theta(u))G_i(u)]$ are continuous functions which attain their minimum and maximum values on this compact set.  By Theorem~\ref{prop monotone}~i), the minimum value of $\frac{\partial M_i}{\partial u_j}$ is positive.  Thus for a large enough value of $Q>0$, $\frac{\partial f_i}{\partial u_j}(u)>0$ on $\R^n$.  Similarly for all other choices of $i,j$, $i\not=j$.

Let $u(t)$ be a solution of (\ref{eq extension}), and let $v(t):=S(u(t))/n$ as in Theorem~\ref{prop monotone}.  Then
\[
nv'(t)=Q\sum_{i=1}^n M_i(u(t))  \ +  \ (1-\theta(u(t)))\sum_{i=1}^n G_i(u(t))
=Q\sum_{i=1}^n M_i(u(t))
\]
\[
=Qn\gamma(S(u(t)))=Qn\gamma(nv(t)),
\]
using the fact that $S(G(u))=0$ on $\R^n$ and (\ref{eq M def}).
Thus it follows that $v'(t)=Q\gamma(n v(t))$, as in (\ref{eq v}).  In particular once more, if $S(u(0))>0$, then $S(u(t))/n\to P$ as $t\to \infty$.  Since $f=Q M$ whenever $S(u)>1/2$, it follows $u(t)\to (P,\ldots, P)$ by Theorem~\ref{prop monotone}, iii).  Similarly for $S(u(0))<0$.

In order to address the fourth statement, we show first that $M(u)\circ u\leq 0$, for all $u\in H$.  To see this, note that for $u\in H$ it holds $S(u)=0$ and thus $M_i(u)=-\theta(u)u_i$, $i=1,\ldots, n$.  Thus $M(u)\circ u=-\theta(u) u\circ u= -\theta(u) \norma{u}^2\leq 0$.

The fact that $f$ has no zeros outside $H$ other than $\pm(P,\ldots, P)$ follows from the argument above.  Suppose that $u\circ g(u)<0$ on $H-R$, and consider $u\in H$, $f(u)=0$.  If $M(u)=0$, then $\theta(u)=0$ by definition of $\theta$ and necessarily $G(u)=g(u)=0$.  If $M(u)\not=0$, then also $\theta(u)<1$, else $0=f(u)=M(u)\not=0$.  In that case also $u\not\in R$ by construction of $M$, hence $u\circ g(u)<0$.  Therefore

\[
f(u)\circ u = (QM(u) + (1-\theta(u)) g(u))\circ u \leq (1-\theta(u)) g(u)\circ u <0,
\]
a contradiction.
\epr






\section{Continuum of Equilibria}  \label{section continuum}

In this section we construct a (non-cooperative) two dimensional reaction diffusion equation with a continuum of spatially nonhomogeneous equilibria, and whose corresponding reaction system (\ref{ode}) is globally attractive towards a single equilibrium.  We begin in reverse by defining the functions which will constitute the nonhomogeneous equilibria, and we build the reaction function based on them.  Consider $\Omega=(-\pi/2,\pi/2)$ and the function

\[
\phi_\lambda(x):=\lambda(  \cos (\sin x),
 \sin (\sin x)), \ x\in \Omega,
\]
defined on $\overline{\Omega}$ for every $\lambda>0$.  In particular, the image of $\phi_\lambda$ is an arc of radius $\lambda$ spanning an angle of two radians.  One can easily compute
$\phi'_\lambda(-\pi/2)=\phi'_\lambda(\pi/2)=0$.  Moreover, 

\be{eq phi 1}
(\phi_\lambda''(x))_1= \lambda \sin(x) \sin(\sin(x)) - \lambda \cos^2 x \cos (\sin x),
\ee

\be{eq phi 2}
(\phi_\lambda''(x))_2= - \lambda \sin(x) \cos(\sin(x)) - \lambda \cos^2 x \sin (\sin x).
\ee

The vector $\phi_\lambda''(x)$ can be thought of as the acceleration vector of $\phi_\lambda(x)$ as $x$ grows from $-\pi/2$ to $\pi/2$, of course, and it points towards the inside of the circle of radius $\lambda$, except for $x=\pm \pi/2$ where it is tangential to this circle.


Let
\[
A:=\{r(\cos \theta, \sin\theta)\,|\, \lambda_1\leq r\leq \lambda_2, -1\leq \theta\leq 1\},
\]
for fixed $0< \lambda_1<\lambda_2$ in $\R^+$.  Given $u\in A$, $u=r(\cos \theta, \sin\theta)$, we denote $\lambda(u):=r$, $x(u):=\sin^{-1}(\theta)$.  Note that these functions are well defined and smooth on $A$, and that $\phi_{\lambda(u)}(x(u))=u$.  Define the vector field $\alpha:A\to \R^2$ by
\be{alpha def}
\alpha(u):=- \phi''_{\lambda(u)}(x(u)).
\ee

\begin{Proposition}  \label{prop continuum}
There exists a smooth function $g:\R^2\to \R^2$ such that i) $g=\alpha$ on the set $A$, and ii) all solutions of the system
\be{ode g}
\frac{d u_i}{d t}=g_i(u),\ \ \ i=1,2,
\ee 
converge towards a single equilibrium.  
\end{Proposition}

\bpr

We start with some basic facts about $\alpha(u)$ for $u\in A$.  It follows from equation (\ref{eq phi 1}) that $\alpha_1(u)>0$ for $u_2=0$ (since $x(u)=0$).  It also follows, from equation (\ref{eq phi 2}), that $\alpha_2(u)>0$, $\alpha_2(u)<0$, and $\alpha_2(u)=0$ whenever $u_2>0$, $u_2<0$, and $u_2=0$, respectively (since $x(u)>0,\ x(u)<0,$ and $x(u)=0$ in each case).   

By (\ref{alpha def}) and the definition of $x(u),\lambda(u)$, the function $\alpha$ is smooth on $A$.  Let $\alpha_1$ and $\alpha_2$ be embedded into smooth functions defined on a closed neighborhood $A'$ of $A$, in such a way that both properties in the previous paragraph still hold on $A'$ for the embedding function, which we also denote by $\alpha$.

Let $A''$ be a closed neighborhood of $A'$.  Consider three smooth functions $\rho_1,\rho_2,\rho_3:\R^2\to [0,1]$ forming a partition of unity of $\R^2$, in the sense that
\begin{enumerate}
\item $\rho_1(u)+\rho_2(u)+\rho_(u)=1$ for every $u\in\R^2$, 
\item $\rho_1(u)\equiv 1$ on $A$, $\rho_1(u)\equiv 0$ on $\R^2-\mbox{int\,}(A')$, and $\rho_1(u)\in (0,1)$ otherwise.
\item $\rho_3(u)\equiv 0$ on $A'$, $\rho_3(u)\equiv 1$ on $\R^2-\mbox{int\,}(A'')$, and $\rho_1(u)\in (0,1)$ otherwise.
\end{enumerate}

In particular, $\rho_2=1-\rho_1$ on $A'-A$, $\rho_2=1-\rho_3$ on $A''-A'$.  

Let $e=(e_1,0)$, for a fixed $e_1>0$ such that $e_1>a_1$ for every $a=(a_1,a_2)\in A''$.   Let $g:\R^2\to \R^2$ be defined by 

\[   
g(u):=\rho_1(u)\alpha(u) + \rho_2(u)(1,0) + \rho_3(u)(e-u).
\]

It is clear that $g$ is smooth and that $g\equiv \alpha$ on $A$ by construction.   We use the Poincare-Bendixson theorem to show that every solution of (\ref{ode g}) must converge towards $e$.  First note that $g\equiv e-u$ outside of a bounded subset of $\R^2$, 
therefore a closed circular region centered on $e$ is invariant and attracts all solutions.  We show below that $g$ has only one equilibrium; therefore every solution must converge towards $e$ or towards some periodic orbit.  But any periodic orbit would contain $e$ inside its enclosed area; see for instance \cite{Strogatz}, Section~6.8.  Moreover, $g_2(u_1,0)=0$ for every $u_1\in \R$, since this is true for each of the vector fields $\alpha$, $u\to e-u$ and $u\to (1,0)$.  Therefore the $u_1$ axis $\{(u_1,0)\,|\, u_1\in \R\}$ is an invariant subset of (\ref{ode g}); this proves that no periodic solutions can exist by the principle of nonintersecting orbits \cite{Strogatz}, and that every solution converges towards $e$.  

It remains to show that $g$ has only one equilibrium.  We show this by considering the various subsets of the partition in $\R^2$.  If $u\in A$, it holds $g(u)=\alpha(u)\not=0$   
by the comments in the beginning of the proof.  If $u\in int(A')-A$, then still $\rho_1(u)>0$, and for $u_2\not= 0$,

\[
g_2(u)=\rho_1(u)\alpha_2(u) + \rho_2(u)\cdot 0 = \rho_1(u)\alpha_2(u)\not=0, \]
by the construction of the extension of $\alpha$ to $A'$.  Also, for $u_2=0$,
\[
g_1(u)=\rho_1(u)\alpha_1(u) + \rho_2(u)\cdot 1 > 0.
\]
If $u\in \R^2 - \mbox{int\,}(A')$, such that $u_1<e_1$, then $\rho_1(u)=0$, and 
\[
g_1(u)= \rho_2(u)\cdot 1 + \rho_3(u)(e_1-u_1) >0.
\]
If $u_1\geq e$, then $g(u)=e-u$, and $g(u)\not=0$ unless $u=e$.  
\epr

\begin{Corollary} Let $g$ be as in Proposition~\ref{prop continuum}, and let $\lambda\in [\lambda_1,\lambda_2]$.  Then the function $\phi_\lambda$ is an equilibrium of the system
\be{pde g}
\frac{\partial u_i}{\partial t}= \Delta u_i + g_i(u),\
 \ \ i=1,2,
 \ee
on $C(\Omega,\R^2)$, under Neumann boundary conditions.
\end{Corollary}

\bpr
By construction of $\alpha$ on $A$, it holds for every $\lambda\in [\lambda_1,\lambda_2]$, $x\in \overline{\Omega}$, that $\alpha(\phi_\lambda(x))=-
\phi''_\lambda(x)$.  Thus $0=\phi''_\lambda(x)+\alpha(\phi_\lambda(x))
=\phi''_\lambda(x)+ f(\phi_\lambda(x))$.  It was calculated that
$\phi_\lambda'(-\pi/2)=\phi_\lambda'(-\pi/2)=0$, $\lambda>0$; this satisfies the required boundary conditions.  
 \epr

\paragraph{Proof of Theorem~\ref{teo appl}:}

\bpr

We use the function $g$ from Proposition~\ref{prop continuum} in the context of Theorem~\ref{teo embedding}.   Let $\pi:\R^2\to H\subseteq \R^3$ be a linear, metric preserving bijection.  Define a smooth function $\tilde{g}:H\to H$ by $\tilde{g}:=\pi\circ g\circ \pi^{-1}$.  Then use Theorem~\ref{teo embedding} to embed this function into a strongly cooperative system (\ref{ode}).  

Recall from the proof of Proposition~\ref{prop continuum} that $g=e-u$ outside of a bounded set.  If $\norma{u}>e_1$, $g(u)=u-e$, then $u\circ g(u)=u\circ e - u\circ u=u_1 e_1 - \norma{u}^2<0$.  Since the bijection $\pi$ preserves angles, this condition is also satisfied for $\tilde{g}$ outside of a bounded region $R$.  Thus we can use item 4. of Theorem~\ref{teo embedding} to conclude that $f$ has a unique equilibrium on $H$.  

It follows automatically from Theorem~\ref{teo embedding} that $\partial f_i/\partial u_j>0$ for every $u$ and $i\not=j$.  Every solution of (\ref{ode}) outside of the H converges towards $\pm (P,\ldots P)$ (Theorem~\ref{teo embedding}, item 3.),  and every solution of (\ref{ode}) in $H$ converges towards $\pi(e)$ by Proposition~\ref{prop continuum}.   Thus every solution of (\ref{ode}) converges towards one of only three equilibria.

Let $u(x):=\phi_\lambda(x)$, $\lambda\in [\lambda_1,\lambda_2]$, be any of the nonhomogeneous equilibria of (\ref{pde g}), so that $u''(x)+ g(u(x))=0$ for every $x\in \Omega$.  Let $\tilde{u}(x):=\pi(u(x))$.  By evaluating $\pi$ on both sides of the previous equation we obtain

\[
0=\pi u''(x) + \pi g(u(x))=\tilde{u}''(x) +\tilde{g}(\tilde{u}(x))=\tilde{u}''(x) +f(\tilde{u}(x)).
\]
Thus $\tilde{u}$ is an equilibrium of (\ref{pde}).  After varying $\lambda$ over $[\lambda_1,\lambda_2]$, the third statement follows.   
\epr

\section{Discussion}

One might imagine a strengthening of Theorem~\ref{teo appl}) in which every solution of the reaction system (\ref{ode}) converges towards a \emph{single} equilibrium $e$.  A standard `sandwich' argument shows that this is impossible:  let (\ref{ode}) be strongly cooperative and converge globally towards $e\in \R^n$.  Then all solutions of the reaction diffusion system (\ref{pde}) must converge towards $e$ as well.   To see this, let $u(x,t)$  be any solution of (\ref{pde}), and let $v\leq u(x,0)\leq w$ componentwise, for some $v,w\in \R^n$.  If $v(t),w(t)$ are the solutions of (\ref{ode}) with initial conditions $v,w$ respectively, then $v(t),w(t)$ also form spatially homogeneous solutions of (\ref{pde}), and $v(t)\leq u(x,t)\leq w(t)$ componentwise for all $t$, by cooperativity \cite{Smith:monotone}.  But $v(t), w(t)$ converge towards $e$, therefore so does $u(t,x)$.)

Regarding item 3. from Theorem~\ref{teo appl}, in a paper by Kishimoto and Weinberger \cite{Kishimoto:1985} it is proved that for convex $\Omega$ and given a system (\ref{pde})
such that $\partial f_i/\partial u_j >0$ for $i\not=j$, any nonhomogeneous equilibrium must be linearly unstable.  Moreover, in a recent result by Smith, Hirsch, and the
author, it is shown that under the same hypotheses the generic
bounded solution of (\ref{pde}) converges towards a homogeneous equilibrium \cite{Enciso:Hirsch:Smith:JDDE2006}.
In particular, the set of nonhomogeneous equilibria must also be sparse (in the sense of prevalence; see \cite{Yorke:1992,Christensen:1972}).  This indicates that item 3. in Theorem~\ref{teo appl} cannot be strenghened to a substantially larger set of nonhomogeneous equilibria.

\paragraph{Acknowledgements}  the author would like to thank Hal Smith for many helpful comments and references.  Also to be thanked are Avner Friedman, Yuan Lou and Eduardo Sontag for their help and encouragement.

\end{document}